\documentclass{amsart}

\usepackage{amssymb}
\usepackage{amsmath,amsthm}
\usepackage{datetime}

\newcommand{\loc}{\ensuremath{\text{loc}}}
\newcommand{\R}{\mathbb R}  

\newtheorem{theorem}{Theorem}

\title{Every Lipschitz metric has ${\mathcal C}^{1}$-geodesics}
\author{Roland Steinbauer}
\address{Faculty of Mathematics, University of Vienna, Oskar-Morgenstern-Platz
1, 1090 Vienna, Austria}
\email{roland.steinbauer@univie.ac.at}
\date{\today}

\begin{document}

\begin{abstract}
We observe that the geodesic equation for any semi-Riemannian metric of
regularity
${\mathcal C}^{0,1}$ possesses
${\mathcal C}^{1}$-solutions in the sense of Filippov.
\vskip 1em

{\em PACS 2010:} 04.20.Cv, 
02.40.Kv, 
02.30.Hq 
 
\end{abstract}

\maketitle

This work is motivated by the recently increased interest in low regularity
issues
in semi-Riemannian geometry and in general relativity. In \cite{CG11} a study of
causality theory of  \emph{continuous} Lorentzian metrics has revealed that some
of the
foundational results do no longer hold in this case. In particular, light cones
may have a non-empty interior and the push up principle ceases to hold. 
On the other hand, standard causality theory holds for ${\mathcal
C}^2$-metrics (\cite{C11}) and was widely believed to still hold in
the \emph{regularity class ${\mathcal C}^{1,1}$} (also denoted by
${\mathcal C}^{2-}$, the metric
being ${\mathcal C}^1$ with the first order derivatives being locally
Lipschitz continuous). Indeed ${\mathcal C}^{1,1}$ is the threshold where
one still has local unique solvability of the geodesic equation. But only
in \cite{KSS13} and independently in \cite{M13} it was shown that the
exponential map of a ${\mathcal C}^{1,1}$-metric possesses the maximal
regularity: it is a bi-Lipschitz homeo\-morphism around each point. Finally
in \cite{M13,KSSV13} it has been shown that causality theory for
${\mathcal C}^{1,1}$-metrics indeed does not deviate from the classical theory.
A more explicit issue was taken up in \cite{LSS13}. There it has been
shown that the geodesics in impulsive pp-wave space-times in Rosen-form, where
the metric is \emph{Lipschitz continuous}, actually are ${\mathcal
C}^{1}$-curves. To this
end the geodesic equation was solved in the sense of Carath\'eodory (see e.g.\
\cite[Ch.\ 1]{F88}) which basically means to solve the classically equivalent
integral equation. 
\medskip

In this short note we observe that, more generally, any 
${\mathcal C}^{0,1}$ (i.e. locally Lipschitz continuous) semi-Riemannian metric
has ${\mathcal C}^{1}$-geodesics. 

In the Riemannian case much more can be shown: The existence of geodesics 
(locally minimizing curves, even in the class of absolutely continuous curves,
see \cite{B12}) is guaranteed for continuous metrics (\cite{H04}) and
if the metric is of H\"older regularity ${\mathcal C}^{0,\alpha}$ ($0<\alpha\leq
1$) they are of regularity ${\mathcal C}^{1,\beta}$ where
$\beta=\alpha/(2-\alpha)$ (\cite{LY06}). If $\alpha=1$ then $\beta=1$ and in
this case the minimizing curves satisfy the geodesic 
equation almost everywhere. 

In the semi-Riemannian case we have to resort to the geodesic equation from the 
start. This, however, is an ODE with discontinuous right hand side and it is
most natural to invoke Filippov's solution concept (\cite{F88}) which is widely
used in non-smooth dynamical systems and non-smooth mechanics (e.g.\
\cite{C08}).

To be more precise, let $M$ be a smooth, $n$-dimensional manifold with a 
${\mathcal C}^{0,1}$-semi-Riemannian metric $g$. By Rademacher's theorem  $g$ is
(classically) differentiable almost everywhere with locally bounded
derivatives. Hence the right hand side of the geodesic equation
\begin{equation}\label{geo}
\ddot x^i=-\Gamma^i_{jk} \dot x^j \dot x^k
\end{equation}
is in $L^\infty_\loc$ but generally not better and we invoke the
concept of Filippov solutions (\cite{F88}).
Indeed when considering an autonomous system of ODEs 
\begin{equation}\label{eq:ode}
 \dot x(t)=F(x(t)) 
\end{equation}
with 
discontinuous right hand side
$F\colon\R^n\supseteq U\to\R^n$ we clearly have to go beyond the classical
${\mathcal C}^1$-solution theory. The key idea of Filippov is not to look
at the values of the discontinuous vector field $F$ at a single point but to
consider a set of directions given by the values of $F$ near that point.
This already hints at the fact that this approach is compatible with
approaches based on regularization, see \cite[Sec.\ 3.3]{H08}. 

Technically Filippov's approach is realized by replacing the ODE by a
\emph{differential inclusion}
\begin{equation}\label{eq:di}
 \dot x(t)\in{\mathcal F}[F](x(t)),
\end{equation}
where the \emph{Filippov set-valued map} ${\mathcal F}[F]\colon \R^n\supseteq
U\to {\mathcal K}_0(U)$ (the collection of all nonempty,
closed and convex subsets of $U$) associated with $F$ is defined as
its essential convex hull, i.e.,
\[
 {\mathcal
F}[F](x):=\bigcap_{\delta>0}\bigcap_{\mu(S)=0}
  {\text{co\,}}\Big(F(B(x,\delta)\setminus S)\Big).
\]
Here ${\text{co}}(B)$ denotes the closed convex hull of a set $B$, $B(x,\delta)$
is the closed ball of radius $\delta$ around $x$, and $\mu$ is the Lebesgue
measure. The explicit calculation of a Filippov set-valued map associated with
a given vector field can be non-trivial, however, there exists a calculus
to simplify this task (\cite{PS87}). Also note that a Filippov set-valued map
is multi-valued only at the points of discontinuity of the original vector
field. 

A \emph{Filippov solution} of \eqref{eq:ode} by definition is an
absolutely continuous curve $x\colon I\to U$, defined on some interval $I$, that
satisfies \eqref{eq:di} almost everywhere. Clearly every classical solution is a
Filippov solution but the latter exist under much milder conditions. We will, 
in particular, make use of the following simple existence result (\cite[Thm.\ 8,
p.\ 85]{F88}).
\begin{theorem}\label{thm:ex}
 If $F\in L^\infty_\loc(U;\R^n)$ then for all $x_0\in U$ there exists a
 Filippov solution of \eqref{eq:ode} with initial condition $x(0)=x_0$.
\end{theorem}

With this preparation we may now state our main result.

\begin{theorem}\label{thm:main}
 Let $(M,g)$ be a smooth manifold with a $\mathcal{C}^{0,1}$-semi-Riemannian
 metric. Then there exist Filippov solutions of the geodesic equation which are 
 ${\mathcal C}^1$-curves.
\end{theorem}

Indeed we only have to rewrite \eqref{geo} as a first order system 
\begin{equation}\label{eq:g1o}
  \dot z(t)=F(z(t))
\end{equation}
in the $2n$ variables $z=(x,\dot x)$ where the right hand side takes the form
$F(z)=(\dot x^1,\dots ,\dot
x^n,-\Gamma^1_{jk}(x)\dot x^j\dot x^k,\dots,-\Gamma^n_{jk}(x)\dot x^j\dot
x^k)$.
Clearly $F\in L^\infty_\loc$, hence by Thm.\ \ref{thm:ex} there exists an
absolutely continuous curve $z$ which solves \eqref{eq:g1o} in the sense of
Filippov. But this implies $\dot x$ is (absolutely) continuous and hence 
$x$ is a $\mathcal{C}^1$-Filippov solution of \eqref{geo}.
\medskip

\emph{Uniqueness} of geodesics, however, is a different issue since the
regularity we are dealing with is clearly below the classical uniqueness
conditions. One possibility is to use \emph{one-sided Lipschitz conditions}.
They will
provide one-sided (either to the future or to the past) uniqueness, see
\cite[\S10.1]{F88} (without forcing the Christoffel symbols to be continuous). 
However, for the main applications we
have in mind, i.e., impulsive gravitational waves, the
metric will be smooth off some null-hypersurface. In this case unique Filippov
solutions do exists under much weaker conditions, see \cite[\S10.2]{F88}. 
Future work will be concerned with using these ideas to prove the
$\mathcal{C}^1$-property of geodesics in non-expanding impulsive gravitational
waves in a cosmological background as well as of expanding impulsive
gravitational waves. 

Finally, future work will also be concerned with relating the present approach
to calculus of variations, in particular, the existence of shortest curves in
the Riemannian case.
\bigskip

{\bf Acknowledgements.} The author was supported by FWF project P25326 and OeAD
project WTZ CZ 15/2013.

\end{document}